\documentclass[12pt] {article}
\usepackage{amsmath,amsthm}
\usepackage{amssymb,latexsym}

\title{Valuations on convex functions and convex sets and Monge-Amp\`ere operators.}
\date{}
\author{Semyon Alesker \footnote{Partially supported by ISF grant 865/16.}
\\  { \normalsize Department of Mathematics, Tel Aviv University, Ramat Aviv}
\\  { \normalsize 69978 Tel Aviv, Israel }
\\ {\normalsize e-mail: semyon@post.tau.ac.il}}

\def\RR{\mathbb{R}}
\def\CC{\mathbb{C}}

\def\HH{\mathbb{H}}

\def\OO{\mathbb{O}}

\def\eps{\varepsilon}

\def\Ome{\Omega}
\def\lam{\lambda}

\def\to{\longrightarrow}

\def\qed { Q.E.D. }

\def\inj{\hookrightarrow}

\swapnumbers
\newtheorem{theorem}{Theorem}[section]

\newtheorem{lemma}[theorem]{Lemma}

\newtheorem{claim}[theorem]{Claim}
\theoremstyle{definition}

\newtheorem{example}[theorem]{Example}
\newtheorem{definition}[theorem]{Definition}
\newtheorem{remark}[theorem]{Remark}
\theoremstyle{conjecture}

 \def\ch{{\cal H}} 
 \def\ck{{\cal K}}

\def\pt{\partial}

\numberwithin{equation}{section}
\begin{document}
\maketitle

\begin{abstract}
The notion of a valuation on convex bodies is very classical. The notion of a valuation on a class of functions was recently introduced and studied by M. Ludwig and others.
We study an explicit relation between continuous valuations on convex functions  which are invariant under adding arbitrary linear functionals, and translations
invariant continuous valuations on convex bodies. More precisely, we construct a natural linear map from the former space to the latter and prove that it has
dense image and infinite dimensional kernel. The proof uses the author's irreducibility theorem and few properties of the real Monge-Amp\`ere operators due to A.D. Alexandrov and Z. Blocki. Furthermore we show how to use complex, quaternionic, and
octonionic Monge-Amp\`ere operators to construct more examples of continuous valuations on convex functions in an analogous way.
\end{abstract}

\section{Introduction.}\label{S:introduction} The main result of the paper is Theorem \ref{T:main-intro} below.
Let $V$ be a finite dimensional real vector space. Let $\ck(V)$ denote the family of all convex compact non-empty subsets of $V$.
\begin{definition}\label{val-sets}
A functional $\phi\colon \ck(V)\to \RR$ is called a valuation if it satisfies the following additivity property:
$$\phi(A\cup B)=\phi(A)+\phi(B)-\phi(A\cap B)$$
whenever $A,B,A\cup B\in \ck(V)$.
\end{definition}

Theory of valuations on convex sets is a classical part of convex geometry which (at least partly) originates in M. Dehn's solution of the third Hilbert problem.
For the classical theory we refer to the surveys \cite{mcmullen-schneider}, \cite{mcmullen-survey} and Ch. 6 in the book \cite{schneider-book}.
In the last two decades there was a considerable progress in the theory: new classification results of various classes of valuations have been proven, new structures on valuations
have been discovered, new applications have been obtained.
This progress and the needs of applications, e.g. to integral geometry, have motivated few attempts to understand the natural generality of some of these developments.
Thus some generalizations of the notion of valuation have been suggested.

One of these recent generalizations is the notion of valuation on a class of functions introduced and studied by M. Ludwig and others, see e.g.
\cite{kone-14}, \cite{ludwig-function-11}, \cite{ludwig-sobolev-12},
\cite{ludwig-colesanti-mussnig}, \cite{tsang-10},
\cite{wang-14}. Let us give a definition of valuation on convex functions which is relevant for the paper.

\begin{definition}\label{val-function}
Let $\Ome \subset V$ be an open subset. An $\RR$-valued functional $\Phi$ on the class of all convex functions  $\Ome\to \RR$ is called a valuation if
$$\Phi(\max\{f,g\})=\Phi(f)+\Phi(g)-\Phi(\min\{f,g\})$$
whenever $f,g,\min\{f,g\}$ are convex in $\Ome$.
\end{definition}
While valuations on convex functions do not apparently look like finitely additive measures, they turn out to be closely related to valuations on convex sets in
the sense of Definition \ref{val-sets}. The main result of this paper (see below) provides some information on this relation in the case of continuous valuations.
To formulate the main result we will need few more definitions.

\begin{definition}\label{D:val-sets-cont}
A valuation on convex compact sets $\phi\colon \ck(V)\to \RR$ is called continuous if it is continuous in the Hausdorff metric on $\ck(V)$.
\end{definition}

\begin{definition}\label{D:val-func-cont}
A valuation $\Phi$ on the set of convex functions in $\Ome$ is called continuous if it is continuous in the $C^0$-topology on this set (e.g. in
the topology of uniform convergence on compact subset of $\Ome$).
\end{definition}

Now any continuous valuation $\Phi$ on convex functions on the dual space $V^*$ induces a continuous valuation $\phi$ on convex compact subsets of $V$ as follows:
\begin{eqnarray}\label{E:map-T}
\phi(K):=\Phi(h_K) \mbox{ for any } K\in\ck(V),
\end{eqnarray}
where $h_K\colon V^*\to \RR$ is the supporting functional of the set $K$ defined by $h_K(\xi):=\sup_{x\in K}\xi(x)$. The continuity of $\phi$ follows from the continuity of $\Phi$ and the well known
(and easy) fact that a sequence of convex compact sets $\{K_i\}$ converges to a convex compact set $K$ if and only if $h_{K_i}\to h_K$ in $C^0$-topology.
The valuation property of $\phi$ follows from the easy fact that if $A,B,A\cup B \in \ck(V)$ then
$$h_{A\cup B}=\max\{h_A,h_B\},\, h_{A\cap B}=\min\{h_A,h_B\}.$$

Let us denote by $Val(V)$ the space of continuous translation invariant valuations on $\ck(V)$
($\phi$ is called translation invariant if $\phi(K+x)=\phi(K)$ for any $K\in\ck(V),x\in V$). $Val(V)$ being equipped with the topology of uniform convergence on
compact subsets of $\ck(V)$ is a Fr\'echet space, in fact even a Banach space.

Let us denote by $VConv(V^*)$ the space of continuous valuations on convex functions on $V^*$ which satisfy the following extra condition
\begin{eqnarray}\label{E:val-func-trans}
\Phi(f+l)=\Phi(f)\, \forall \mbox{ convex function } f \mbox{ and } \forall \mbox{ linear functional } l \mbox{ on } V^*.
\end{eqnarray}

The correspondence $\Phi\mapsto \phi$ as in (\ref{E:map-T}) gives readily a linear map
$$T\colon VConv(V^*)\to Val(V).$$
Now we can formulate the main result of the paper.
\begin{theorem}\label{T:main-intro}
The map $T\colon VConv(V^*)\to Val(V)$ has the image dense in $Val(V)$. Its kernel is infinite dimensional.
\end{theorem}
Let us say a few words on the proof of this result. The map $T$ commutes with the natural actions on both spaces of group $GL(V)$ of
invertible linear transformations on $V$.
Hence by the author's irredicibility theorem \cite{alesker-irr} it suffices to show that $Im(T)$ intersects non-trivially each subspace of $Val(V)$ consisting
of valuations of the given degree of homogeneity and parity. In order to construct sufficiently many such examples we use the real Monge-Amp\`ere operator.
For general convex functions it was defined by A.D. Alexandrov \cite{alexandrov2} who also proved the necessary continuity property. The valuation property was
proved by Blocki \cite{blocki} in a somewhat different language. (In fact Blocki originally proved the corresponding property for the complex Monge-Amp\`ere operator, but
the real case is easily implied by it.)

After proving the main result in Section \ref{S:proof}, we discuss in the subsequent section analogous constructions of valuations on convex functions using complex,
quaternionic, and octonionic Monge-Amp\`ere operators (the latter is only in the case of 2 octonionic variables). While the real and complex Monge-Amp\`ere operators
are very classical, the quaternionic and octonionic ones were introduced relatively recently by the author \cite{alesker-bsm}, \cite{alesker-octon}.

\section{Proof of the main result.}\label{S:proof}
 \underline{Step 1.} In this step we show that the image of $T$ is a $GL(V)$-invariant subspace of $Val(V)$.
Observe that the map $T$ commutes with the natural action of the group $GL(V)$ on both spaces.
Let us define these actions explicitly. Let $\Phi\in VConv(V^*),\phi\in Val(V),g\in GL(V)$. Then
\begin{eqnarray*}
(g\Phi)(f)=\Phi(f\circ g^{*-1}) \mbox{ for any convex } f\colon V^*\to \RR,\\
(g\phi)(K)=\phi(g^{-1}K) \mbox{ for any } K\in \ck(V).
\end{eqnarray*}

\hfill

\underline{Step 2.} Remind the McMullen's decomposition \cite{mcmullen-decomp} of $Val(V)$. Let us denote by $Val_i(V)$ the subspace of $Val(V)$ of $i$-homogeneous valuations, i.e.
$$Val_i(V):=\{\phi\in Val(V)|\, \phi(\lam K)=\lam^i\phi(K) \mbox{ for any }\lam>0 \mbox{ and }K\in \ck(V)\}.$$
Then McMullen's theorem \cite{mcmullen-decomp} says that
\begin{eqnarray}\label{E:homog}
Val(V)=\oplus_{i=0}^{\dim V} Val_i(V).
\end{eqnarray}
Furthermore $Val_i(V)$ can be decomposed further with respect to parity, i.e.
\begin{eqnarray}\label{E:parity}
Val_i(V)=Val_i^+(V)\oplus Val_i^-(V),
\end{eqnarray}
where $Val_i^+(V)$ is the subspace of even valuations, i.e. satisfying $\phi(-K)=\phi(K)$ for any $K\in \ck(V)$, and $Val_i^-(V)$
is the subspace of odd valuations, i.e. satisfying $\phi(-K)=-\phi(K)$ for any $K\in \ck(V)$. The group $GL(V)$ preserves degree of homogeneity and parity of valuations
from $Val(V)$, i.e. the two decompositions (\ref{E:homog})-(\ref{E:parity}). Furthermore by the author's irreducibility theorem \cite{alesker-irr} the action of $GL(V)$ on each
$Val_i^\pm(V)$ is topologically irreducible, i.e. any non-zero $GL(V)$-invariant subspace is dense.

Hence it suffices to show the following claim.
\begin{claim}\label{Cl:main-claim}
$Im(T)\cap Val_i^\pm(V)\ne 0$ whenever $Val_i^\pm(V)\ne 0$.
\end{claim}

\underline{Step 3.} Before we prove Claim \ref{Cl:main-claim} we have to remind few basic notions. First remind the notion of mixed determinant.

Let $W$ be a finite dimensional real vector space, and let $p\colon W\to \RR$ be a homogeneous polynomial of degree $n$. Then it is well known (and easy to see)
that there exists a unique $n$-linear symmetric map
$$\bar p\colon \underset{n \mbox{ times}}{\underbrace{W\times\dots\times W}}\to \RR$$
such that $\bar p(x,\dots,x)=p(x)$ for any $x\in W$.

Let us apply this construction to the space $W=\ch_n(\RR)$ of real symmetric matrices of size $n$ and $p=\det$ be the determinant. The corresponding
$n$-linear symmetric map is called the mixed determinant and is denoted by $\det$ (without the bar).

\hfill

\underline{Step 4.} Let us remind a couple of facts on the real Monge-Amp\`ere operators. For any $C^2$-smooth function $f$ on $\RR^n$ we denote by
$Hess_\RR(f)=\left(\frac{\pt^2f}{\pt x_i\pt x_j}\right)_{i,j=1}^n$ its real Hessian.

\begin{theorem}[Alexandrov \cite{alexandrov2}]\label{T:alexandrov}
Let $\Ome\subset \RR^n$ be an open subset. Let $1\leq i\leq n$. Let $A_1,\dots,A_{n-i}$ be continuous compactly supported
functions on $\Ome$ with values in the space $\ch_n(\RR)$ of real symmetric matrices of size $n$. Then for any convex function $f\colon \Ome\to \RR$ one can define a
signed measure denoted by
\begin{eqnarray*}
\det(\underset{i\mbox{ times}}{\underbrace{Hess_\RR(f),\dots,Hess_\RR(f)}},A_1(x),\dots,A_{n-i}(x))\equiv\\ \det(Hess_\RR(f)[i],A_1(x),\dots,A_{n-i}(x))
\end{eqnarray*}
which satisfies the following properties:

(i) if $f$ is $C^2$-smooth then this measure has the obvious meaning of the mixed determinant of matrix-valued functions.

(ii) if a sequence $\{f_l\}_{l=1}^\infty$ of convex in $\Ome$ functions converges to $f$ in $C^0$-topology, i.e. uniformly on compact subsets of $\Ome$,
then $$\det(Hess_\RR(f_l)[i],A_1(x),\dots,A_{n-i}(x))$$ converges weakly in the sense of measures to
$$\det(Hess_\RR(f)[i],A_1(x),\dots,A_{n-i}(x)).$$

(iii) For any open subset $U\subset \Ome$ the construction of the measure commutes with the restriction to $U$, i.e. one has equality of measures on $U$
\begin{eqnarray*}
\det(Hess_\RR(f|_U)[i],(A_1|_U)(x),\dots,(A_{n-i}|_U)(x))=\\
\left(\det(Hess_\RR(f)[i],A_1(x),\dots,A_{n-i}(x))\right)\big|_U.
\end{eqnarray*}

(iv) If the functions $A_1,\dots,A_{n-i}$ are non-negative definite matrices pointwise then the measure $\det(Hess_\RR(f)[i],A_1(x),\dots,A_{n-i}(x))$ is non-negative.
\end{theorem}

We will need another result essentially due to Blocki \cite{blocki}. Strictly speaking, Blocki proved it for the complex Monge-Amp\`ere operator and in a slightly different (in fact, more general) form.
However the real case follows immediately from the complex one by imbedding $\RR^n\inj \CC^n=\RR^n\oplus \sqrt{-1}\RR^n$. The needed form of the result is deduced in \cite{alesker-aim05} from the Blocki's form
of it in a slightly more general situation of quaternionic Monge-Amp\`ere operator.
\begin{theorem}[Blocki \cite{blocki}]\label{T:blocki}
Let $\Ome\subset \RR^n$ be an open subset. Let $1\leq i\leq n$. Let $A_1,\dots,A_{n-i}$ be continuous compactly supported
functions on $\Ome$ with values in the space $\ch_n(\RR)$ of real symmetric matrices of size $n$. Let $f,g\colon\Ome\to \RR$ be convex functions
such that $\min\{f,g\}$ is also convex. Then the measures from Theorem \ref{T:alexandrov} satisfy the following valuation property:
\begin{eqnarray*}
\det(Hess_\RR(\min\{f,g\})[i], A_1(x),\dots,A_{n-i}(x))+\\\det(Hess_\RR(\max\{f,g\})[i], A_1(x),\dots,A_{n-i}(x))=\\
\det(Hess_\RR(f)[i], A_1(x),\dots,A_{n-i}(x))+\det(Hess_\RR(g)[i], A_1(x),\dots,A_{n-i}(x)).
\end{eqnarray*}
\end{theorem}

Theorems \ref{T:alexandrov} and \ref{T:blocki} immediately imply that for any $1\leq i\leq n$ for any continuous compactly supported functions
$A_1,\dots,A_{n-i}\colon \Ome\to \ch_n(\RR)$ and $B\colon \Ome \to \RR$ the functional
\begin{eqnarray}\label{E:real-val-f}
\Phi(f):=\int_\Ome B(x)\det(Hess_\RR(f)[i], A_1(x),\dots,A_{n-i}(x))dx
\end{eqnarray}
is a valuation from $VConv(\RR^n)$.

The corresponding valuation $\phi\in Val(\RR^n)$ is
\begin{eqnarray}\label{E:homog-val}
\phi(K)=\int_\Ome B(x)\det(Hess_\RR(h_K)(x)[i], A_1(x),\dots,A_{n-i}(x))dx
\end{eqnarray}
for any $K\in \ck(\RR^n)$. Obviously $\phi$ is $i$-homogeneous.

\hfill

\underline{Step 5.} Recall that $Val_0(V)$ is 1-dimensional and is spanned by the Euler characteristic $\chi$ (this statement is trivial).
Clearly $\chi=T(\tilde\chi)$ where $\tilde\chi\in VConv(V)$ is defined by $\tilde\chi(f)=1$ for any convex function $f$. Hence $Val_0(V)\subset Im(T)$.

\hfill

\underline{Step 6.} Let $1\leq i\leq n-1$. In this case $Val_i^\pm(V)$ is known to be infinite dimensional.
Let us show that $Im(T)\cap Val_i^\pm(V)\ne 0$ in this case. It suffices to show that there exist non-zero valuations of the form (\ref{E:homog-val})
of any parity. Since $Im(T)$ is clearly invariant under the involution which fixes even valuations $Val^+(V)$ and changes the sign of odd valuations $Val^-(V)$,
it suffices to construct a non-zero valuation of the form (\ref{E:homog-val}) which is neither even nor odd.
This will be achieved by choosing $A_1,\dots,A_{n-i}$ appropriately and $B\equiv 1$.

For $p=1,\dots,n$ let us denote
$$E_p=diag(0,\dots 0,1,0,\dots,0),$$
the diagonal matrix of size $n\times n$ where 1 is located at the $p$th place. Let us denote $v_0=(1,0\dots,0)$.
Let us fix a continuous compactly supported function $\psi\colon \RR^n\to \RR$ such that $\psi(v_0)=1$.

Let us denote
\begin{eqnarray*}
A_1=E_1\delta_{v_0},\\
A_l=E_{l}\psi \mbox{ for } 2\leq l\leq n-i,
\end{eqnarray*}
where $\delta_{v_0}$ is the delta-function on $\RR^n$ supported at $v_0$. We will construct a
convex set $K\in \ck(\RR^n)$ with the supporting functional $h_K$ infinitely smooth outside of 0 such that
\begin{eqnarray*}
\phi(K):=\int_{\RR^n}\det(Hess_\RR(h_K)(x)[i], A_1(x),\dots,A_{n-i}(x))dx=\\\det(Hess_\RR(h_K(v_0))[i],E_1,\dots,E_{n-i})
\end{eqnarray*}
satisfies $\phi(-K)\ne \pm \phi(K)$. Then approximating the delta-function by continuous functions we will get a valuation $\phi'$ of the form
(\ref{E:homog-val}) with continuous compactly supported $A_l$'s such that $\phi(-K)\ne \pm \phi'(K)$, i.e. $\phi'$ is neither even nor odd as required.

Next observe that for any $n\times n$ matrix $H$ one has
$$\det(H[i],E_1,\dots,E_{n-i})={{n}\choose{i}}^{-1}\det H_{i},$$
where by $H_i$ we denote the matrix of size $i\times i$ which is obtained from $H$ by deleting the first $n-i$ rows and columns.

Hence we get that
\begin{eqnarray}\label{E:001}
\phi(K)={{n}\choose{i}}^{-1} \det(Hess_\RR (h_K(v_0)))_i.
\end{eqnarray}

Observe that for any 1-homogeneous smooth function $h\colon \RR^n\backslash\{0\}\to \RR$ the vector $v_0$ belongs to the kernel of the matrix
$Hess_\RR(h)(v_0)$. Hence the first row and column of $Hess_\RR(h_k(v_0))$ vanish, but they are deleted anyway in the right hand side of (\ref{E:001}).

Since $\phi(-K)={{n}\choose{i}}^{-1} \det(Hess_\RR (h_K(-v_0)))_i$ it suffices to choose a convex set $K\in \ck(\RR^n)$ with infinitely smooth $h_K$ outside of 0
such that
\begin{eqnarray}\label{E:002}
\det(Hess_\RR (h_K(v_0)))_i\ne \pm \det(Hess_\RR (h_K(-v_0)))_i.
\end{eqnarray}
Such  $K$ can be constructed by first taking the intersection of two Euclidean balls of radii 1 and 2 respectively whose centers are close enough and are located on the line
spanned by $v_0$, and then smoothing the obtained body near the intersection of the boundaries of the balls. For the obtained body $K$ one has
\begin{eqnarray*}
Hess_\RR (h_K(v_0))=diag(0,\underset{n-1 \mbox{ times}}{\underbrace{1,\dots,1}}),\\
Hess_\RR (h_K(-v_0))=2 \cdot diag(0,\underset{n-1 \mbox{ times}}{\underbrace{1,\dots,1}})
\end{eqnarray*}
Hence
$$\det(Hess_\RR (h_K(v_0)))_i=1,\, \det(Hess_\RR (h_K(-v_0)))_i=2^i,$$
hence (\ref{E:002}) is satisfied. Thus Theorem \ref{T:main-intro} is proved for $1\leq i\leq n-1$.

\hfill

\underline{Step 7.} Let us consider the case of the degree of homogeneity $i=n$. By a theorem of Hadwiger \cite{hadwiger-vol}
$Val_n(\RR^n)$ is spanned by the Lebesgue measure $vol_n$. In particular $Val_n^-(\RR^n)=0$.

Let us consider valuations on convex functions on $\RR^n$ of the form
$$\Phi(f)=\int_{\RR^n}B(x)\det Hess_\RR(f) dvol(x),$$
where $B\colon \RR^n\to \RR$ is a continuous compactly supported function.

\begin{lemma}\label{L:003}
For any convex set $K\in \ck(\RR^n)$ one has
$$\det h_K= vol_n(K)\delta_0, $$
where $\delta_0$ is the delta-measure supported at the origin.
\end{lemma}
We will leave the details of the proof of this lemma to the reader, while we will comment on it anyway.
First we may assume by approximation that $h_K$ is smooth outside of 0. Since $h_K$ is 1-homogeneous, for any $x\in \RR^n\backslash\{0\}$
the kernel of $Hess_\RR(h_K(x))$ contains $x$. Hence the measure $\det Hess_\RR (h_K)$ vanishes on $\RR^n\backslash\{0\}$. Hence it must be proportional to $\delta_0$:
$$\det Hess_\RR h_K=C(K)\delta_0.$$
The functional $K\mapsto C(K)$ is continuous in the Hausdorff metric on $\ck(\RR^n)$ by Theorem \ref{T:alexandrov}, and it is a valuation by Theorem \ref{T:blocki}.
Obviously it is also $n$-homogeneous. Hence by the Hadwiger theorem \cite{hadwiger-vol} $C(K)$ is proportional to $vol_n(K)$ with obviously a non-negative constant
of proportionality. As a more direct argument shows, this constant of proportionality is equal to 1.

Lemma \ref{L:003} implies that
\begin{eqnarray}\label{E:004}
\Phi(h_K)=B(0)\cdot vol_n(K).
\end{eqnarray}

Choosing $B$ such that $B(0)\ne 0$ we get the result.

\hfill

\underline{Step 8.} It remains to show that $Ker(T)$ is infinite dimensional. Let us consider an arbitrary continuous compactly
supported function $B\colon \RR^n\to \RR$ such that $B(0)=0$.  The valuation on convex functions on $\RR^n$
$$\Phi(f):=\int_{\RR^n}B(x)\cdot \det (Hess_\RR (f))dvol(x)$$
belongs to $Ker(T)$ by (\ref{E:004}). Let us show that $\Phi= 0$ iff $B\equiv 0$.

Let us assume $\Phi=0$, i.e. $\int_{\RR^n}B(x)\cdot \det (Hess_\RR f(x))dvol(x)=0$ for any convex function $f\colon\RR^n\to \RR$. Let us take $f$ to be of the form
$$f(x)=\frac{|x|^2}{2}+\eps\psi(x),$$
where $\psi\colon \RR^n\to\RR$ is a compactly supported infinitely smooth function and $\eps$ is small enough so that $f$ is convex.
Then $Hess_\RR(f(x))=I_n+\eps \cdot Hess_\RR(\psi(x))$. Hence
$$\det(Hess_\RR(f(x)))=1+\eps\cdot \Delta \psi(x)+O(\eps^2).$$
Hence we have
$$0=\Phi(f)=\int_{\RR^n}(1+\eps\cdot \Delta \psi+O(\eps^2))\cdot B(x) dvol(x).$$
Hence $$\int_{\RR^n}\Delta \psi \cdot B(x) dvol(x)=0$$
for any compactly supported smooth function $\psi$. Integrating by parts we obtain that $\Delta B=0$.
But since $B$ is compactly supported we get $B\equiv 0$.
\qed

\section{Related constructions of continuous valuations on convex functions.}\label{S:further-constructions}
In the previous section we have used a couple of non-trivial properties of the real Monge-Amp\`ere operator to construct continuous valuations on convex functions.
In this section we will briefly outline the use of the other classical complex Monge-Amp\`ere operator and more recently introduced by the author quaternionic and octonionic
Monge-Amp\`ere operators to analogous, though different, constructions of continuous valuations on convex functions.

\subsection{Complex Monge-Amp\`ere operator.}\label{Ss:complex} Let $\Ome\subset \CC^n$ be an open subset.
For a $C^2$-smooth function $f\colon \Ome\to \RR$ the complex Hessian
is defined by
$$Hess_\CC(f(x)):=\left(\frac{\pt^2 f}{\pt z_i\pt\bar z_j}\right),$$
where $z_p=x_p+\sqrt{-1}y_p$ with $x_p,y_p\in\RR$, and
\begin{eqnarray*}
\frac{\pt}{\pt \bar z_i}=\frac{1}{2}\left(\frac{\pt}{\pt x_i}+\sqrt{-1}\frac{\pt}{\pt y_i}\right),\\
\frac{\pt}{\pt z_j}=\frac{1}{2}\left(\frac{\pt}{\pt x_j}-\sqrt{-1}\frac{\pt}{\pt y_j}\right).
\end{eqnarray*}

Let $1\leq i\leq n$. Let $A_1,\dots,A_{n-i}\colon \Ome\to \ch_n(\CC)$ be continuous compactly supported functions with values in the space of complex Hermitian $n\times n$ matrices.
Let $B\colon \Ome\to \RR$ be continuous and compactly supported.
For a $C^2$-smooth function $f\colon\Ome\to \RR$ let us denote by
$$\det\left(Hess_\CC(f(x))[i],A_1(x),\dots,A_{n-i}(x)\right)$$
the obvious mixed determinant. As a generalization of the Alexandrov's theorem \ref{T:alexandrov}, it was shown by Chern, Levine, and Nirenberg \cite{chern-levine-nirenberg}
that this expression can be defined, as a signed measure, for any convex (more generally, continuous plurisubharmonic) function $f$ and this measure
satisfies all the properties stated in Theorem \ref{T:alexandrov}, where $Hess_\RR$ should be replaced with $Hess_\CC$ everywhere.

Furthermore Blocki \cite{blocki} has shown that
\begin{eqnarray}\label{E:com-val}
f\mapsto \int_{\CC^n}B(x) \cdot \det\left(Hess_\CC(f(x))[i],A_1(x),\dots,A_{n-i}(x)\right)dvol(x)
\end{eqnarray}
is a valuation on convex functions (in fact even on continuous plurisubharmonic functions on $\CC^n$).
Thus we get that (\ref{E:com-val}) is a continuous valuation on convex functions on $\CC^n$. Clearly it is invariant under addition of linear functionals.
Thus altogether we get an element of $VConv(\CC^n)$.

\begin{remark}\label{kazarn-complex}
Valuations of the form (\ref{E:com-val}) being restricted to supporting functions of convex polytopes were considered by Kazarnovskii \cite{kazarnovskii-81}, \cite{kazarnovskii-84} (using a different notation)
in the context of complex analysis rather than valuations theory. He was interested in asymptotics of a number of zeros of a system of exponential sums in $n$ variables in terms
of Newton polytopes of the sums. If one takes $i=n$ (thus there is no $A$'s) and the function $B$ to be rotation invariant then all such valuations on convex bodies are proportional
to so called Kazarnovskii's pseudovolume. This pseudovolume is a continuous translation invariant $U(n)$-invariant valuation on convex bodies. The recent progress in understanding of the structure of
such valuations allowed Bernig and Fu to give an alternative description of it in integral geometric terms (see \cite{bernig-fu}, Lemma 3.3 and a remark after it).
\end{remark}

\subsection{Quaternionic  Monge-Amp\`ere operator.}\label{Ss:quatern} Quaternionic Monge-Amp\`ere operator was first defined by the author in \cite{alesker-bsm}.
Let $\HH$ denotes the (non-commutative) field of quaternions.  Any quaternion $q\in \HH$ can be uniquely written in the form
$$q=t+ix+jy+kz,$$
where $t,x,y,z\in\RR$, and $i,j,k$ are the usual anti-commuting quaternionic units satisfying
$$i^2=j^2=k^2=-1, ij=k.$$
(All other standard relations follow from these ones and anti-commutativity.)

\hfill

In order to construct valuations on convex functions on $\HH^n$ similar to the real and complex cases (\ref{E:real-val-f}) and (\ref{E:com-val})
we will define, following \cite{alesker-bsm}, quaternionic Hessian and take mixed determinant of quaternionic Hermitian matrices when determinant is
understood in the sense of Moore (to be described).

Let $\Ome\subset \HH^n$ be an open subset. For any smooth function $F\colon \Ome\to \HH$ define
\begin{eqnarray}\label{E:quat-CR}
\frac{\pt F}{\pt \bar q_a}=\frac{\pt F}{\pt t_a}+i\frac{\pt F}{\pt x_a}+j\frac{\pt F}{\pt y_a}+k\frac{\pt F}{\pt z_a},\\
\frac{\pt F}{\pt  q_b}=\frac{\pt F}{\pt t_b}-\frac{\pt F}{\pt x_b}i-\frac{\pt F}{\pt y_b}j-\frac{\pt F}{\pt z_b}k.
\end{eqnarray}
Note the change of signs and order of terms in the second row. It is straightforward to check using associativity of the product of quaternions that these operators commute:
$$\left[\frac{\pt }{\pt \bar q_a},\frac{\pt }{\pt  q_b}\right]=0.$$

Now let us define the quaternionic Hessian of a $C^2$-smooth function $f\colon \Ome\to \RR$ by
$$Hess_\HH (f(x)):=\left(\frac{\pt^2f}{\pt \bar q_i\pt q_j}\right).$$
This matrix takes quaternionic values and since $f$ is real valued it is Hermitian. Recall that a quaternionic matrix $(a_{ij})_{n\times n}$
is called Hermitian if $a_{ji}=\overline{a_{ij}}$ for any $i,j$, where $\bar a$ denotes the quaternionic conjugation of $a\in \HH$. Let us denote by $\ch_n(\HH)$ the space of
quaternionic Hermitian matrices of size $n$.

Over non-commutative fields there is no notion of determinant of matrices which would have all the properties of the usual determinant in the
commutative case. However there is a notion of the Dieudonn\'e determinant which in the case of quaternions behaves like the absolute
value of the usual real or complex determinant (see Section 1.2 in \cite{alesker-bsm} and references therein). More importantly for this paper, on quaternionic Hermitian
matrices there is a notion of the Moore determinant which has many of the properties of the usual determinant on the real symmetric and complex hermitian matrices.
For example in terms of this determinant one can formulate and prove the Sylvester criterion of positive definiteness of quaternionic Hermitian matrices and
Alexandrov's inequalities for mixed determinants, see Section 1.1 in \cite{alesker-bsm} and references therein.

\hfill

For any $n\times n$ quaternionic matrix $A$ let us define its realization ${}^{\RR} A$ which is a real matrix of size $4n\times 4n$.
Consider the $\RR$-linear operator $\hat A\colon \HH^n\to\HH^n$ defined by multiplication by $A$, i.e. $\hat A(x)=Ax$. If we identify in the standard way
$\HH^n\simeq \RR^{4n}$ we get an $\RR$-linear operator $\tilde A\colon \RR^{4n}\to \RR^{4n}$. Its matrix in the standard basis is called the realization of $A$ and is denoted
by ${}^{\RR} A$. For the following result we refer to \cite{aslaksen} and references therein.

\begin{theorem}
There exists a polynomial $P\colon \ch_n(\HH)\to \RR$ which is uniquely characterized by the following two properties:

(i) $\det({}^{\RR} A)=P(A)^4$ for any $A\in \ch_n(\HH)$, where $\det $ on the left denotes the usual determinant of real matrices of size $4n$;

(ii) $P(I_n)=1$.
\end{theorem}
Obviously $P$ is homogeneous of degree $n$. This polynomial $P$ is called Moore determinant and is denoted by $\det$.
There is no abuse of notation due to the following examples and properties of the Moore determinant for the detailed discussion of which we
refer to \cite{alesker-bsm} and \cite{aslaksen}.
\begin{example}
(1) Any \itshape{complex} Hermitian matrix can be considered as quaternionic Hermitian. Then its Moore determinant equals to its usual determinant of
complex matrices.

(2) General quaternionic Hermitian matrix of size 2 has the form
$$A=\left[\begin{array}{cc}
            a&q\\
            \bar q&b
          \end{array}\right],$$
          where $a,b\in\RR,q\in \HH$. The its Moore determinant $\det A=ab-\bar q q(=ab-q\bar q)$.
(3) Let $A$ be quaternionic Hermitian matrix, and $C$ be any quaternionic matrix of the same size. Then the Moore determinant satisfies the following property
of weak multiplicativity:
$$\det(C^*AC)=\det A\cdot \det(C^*C),$$
where $C^*$ denotes the conjugate matrix obtained from $C$ by transposing and quaternionic conjugation of all elements; observe that $C^*AC$ and $C^*C$ are Hermitian.
\end{example}

Now we can proceed similarly to the real and complex cases. Since the Moore determinant is an $n$-homogeneous polynomial, one can consider mixed Moore determinant.
Let $1\leq i\leq n$. Let $A_1,\dots,A_{n-i}\colon \Ome\to \ch_n(\HH)$ be continuous compactly supported functions with values in the space of quaternionic Hermitian $n\times n$ matrices.
Let $B\colon \Ome\to \RR$ be continuous and compactly supported.
For a $C^2$-smooth function $f\colon\Ome\to \RR$ let us denote by
$$\det\left(Hess_\HH(f(x))[i],A_1(x),\dots,A_{n-i}(x)\right)$$
the obvious mixed determinant. As a generalization of Theorems of Alexandrov \ref{T:alexandrov} and Chern-Levine-Nirenberg mentioned in Section \ref{Ss:complex}, it was shown by the author \cite{alesker-bsm}
that this expression can be defined, as a signed measure, for any convex (more generally, continuous quaternionic plurisubharmonic) function $f$ and this measure
satisfies all the properties stated in Theorem \ref{T:alexandrov}, where $Hess_\RR$ should be replaced with $Hess_\HH$ everywhere.

Furthermore the author \cite{alesker-aim05} has generalized Blocki's Theorem \ref{T:blocki} to quaternionic situation
\begin{eqnarray}\label{E:quat-val}
f\mapsto \int_{\HH^n}B(x) \cdot \det\left(Hess_\HH(f(x))[i],A_1(x),\dots,A_{n-i}(x)\right)dvol(x)
\end{eqnarray}
is a valuation on convex functions on $\HH^n$ (in fact even on continuous quaternionic plurisubharmonic functions on $\HH^n$).
Thus we get that (\ref{E:quat-val}) is a continuous valuation on convex functions on $\HH^n$. Clearly it is invariant under addition of linear functionals.
Thus altogether we get an element of $VConv(\HH^n)$.

\begin{remark}\label{kazarn-quat}
If one takes $i=n$ in (\ref{E:quat-val}) (thus there is no $A$'s) and the function $B$ to be rotation invariant and restricts these expressions
to supporting functions of convex compact sets, then all such valuations on convex bodies are proportional to each other and are considered to be a quaternionic version
of Kazarnovskii's pseudovolume (see Remark \ref{kazarn-complex}).
\end{remark}

\subsection{Octonionic Monge-Amp\`ere operator.}\label{Ss:octon}
The octonionic Monge-Amp\`ere operator was  defined by the author \cite{alesker-octon} for functions of two octonionic variables; it is discussed below.
The octonionic Hessian of a smooth real valued function can apparently be defined for any number of variables in analogy to the real, complex, and
quaternionic situations discussed above; it takes values in octonionic Hermitian matrices. However the next step of taking the determinant of it seems to be problematic in general.
Octonionic Hermitian matrices
of size 2 and probably 3 do admit a good notion of determinant (the case of size 2 is discussed below), but the author is not aware of a good octonionic determinant
in higher dimensions.

\hfill

The non-commutative and non-associative field of octonions will be denoted by $\OO$. Any octonion $q\in\OO$ can be written in the standard form
$$q=\sum_{i=0}^7x_i e_i,$$
where $x_i\in \RR$ and $e_0=1,e_1,\dots, e_7$ are the standard octonionic units (see \cite{alesker-octon} and references therein).
Let $\Ome\subset \OO^2$ be an open subset. If $(q_1,q_2)$ are octonionic coordinates in $\OO^2$ we write $q_a=\sum_{i=0}^7x_{ai} e_i,$ for $a=1,2$.
For any smooth function $F\colon \Ome\to \OO$ let us denote for $a,b=1,2$
\begin{eqnarray*}
\frac{\pt F}{\pt \bar q_a}=\sum_{i=0}^7e_i\frac{\pt F}{\pt x_{ai}},\\
\frac{\pt F}{\pt  q_b}=\sum_{i=0}^7\frac{\pt F}{\pt x_{bi}}\bar e_i,
\end{eqnarray*}
where $\bar q$ is the octonionic conjugation, namely it is the only $\RR$-linear operation such that
$$\bar e_0=e_0, \bar e_i=-e_i \mbox{ for } i>0.$$
It is easy to check that if $f\colon \Ome\to \RR$ is real valued and $C^2$-smooth then
\begin{eqnarray*}
\frac{\pt }{\pt \bar q_a}\frac{\pt f}{\pt  q_b}=\frac{\pt }{\pt  q_b}\frac{\pt f}{\pt \bar q_a} \mbox{ for any } a,b=1,2.
\end{eqnarray*}
For a real valued $C^2$-smooth function $f$ one defines its octonionic Hessian to be
$$Hess_{\OO}(f)=\left(\frac{\pt ^2f}{\pt \bar q_i \pt  q_j}\right)_{i,j=1}^2$$
which is an octonionic Hermitian matrix.

\hfill

Let $\ch_2(\OO)$ denote the space of octonionic Hermitian matrices of size 2. Observe that $\ch_2(\OO)$ consists of elements of the form
\begin{eqnarray*}
\left[\begin{array}{cc}
          a&q\\
          \bar q&b
          \end{array}\right] \mbox{ where } a,b\in \RR, q\in \OO.
\end{eqnarray*}
The (Moore) determinant $\det \colon\ch_2(\OO)\to \RR$ is defined by
$$\det\left[\begin{array}{cc}
          a&q\\
          \bar q&b
          \end{array}\right]=ab-q\bar q(=ab-\bar q q).$$
Since $\det$ is a homogeneous polynomial of degree 2, one can define the mixed determinant.
This $\det$ also satisfies the Sylvester criterion of positive definiteness and the Alexandrov's inequality for mixed determinants \cite{alesker-octon}.

\hfill

Let $A\colon \Ome\to \ch_2(\OO)$ be a continuous compactly supported function. The author \cite{alesker-octon} has generalized the Alexandrov Theorem \ref{T:alexandrov} and the
Chern-Levin-Nirenberg Theorem mentioned in Section \ref{Ss:complex} to the case of two octonionic variables. More precisely for any convex function $f\colon \Ome\to \RR$  one can define the signed measures
$$\det(Hess_\OO(f),A(x)),\, \det(Hess_\OO(f))$$
which satisfy all the properties from Theorem \ref{T:alexandrov}, where $Hess_\RR$ should be replaced with $Hess_\OO$.

Furthermore the author \cite{alesker-octon} has generalized Blocki's Theorem \ref{T:blocki} to octonionic situation, namely
\begin{eqnarray}\label{E:oct-val1}
f\mapsto \int_{\OO^2}B(x) \cdot \det\left(Hess_\OO(f),A(x))\right)dvol(x)\\\label{E:oct-val2} \mbox{ and } f\mapsto \int_{\OO^2}B(x)\cdot \det (Hess_\OO(f)) dvol(x)
\end{eqnarray}
are valuations on convex functions (in fact even on continuous octonionic plurisubharmonic functions on $\OO^2$).
Thus we get that (\ref{E:oct-val1}) and (\ref{E:oct-val2}) are continuous valuations on convex functions on $\OO^2$. Clearly they are invariant under addition of linear functionals.
Thus altogether we get that (\ref{E:oct-val1}) and (\ref{E:oct-val2}) are elements of $VConv(\OO^2)$.
\begin{remark}\label{kazarn-oct}
If one takes in (\ref{E:oct-val2}) the functions $B$ to be rotation invariant and restricts these expressions
to supporting functions of convex compact sets, then all such valuations on convex bodies are proportional to each other and are considered to be an octonionic version
of Kazarnovskii's pseudovolume (see Remark \ref{kazarn-complex}).
\end{remark}

\end{document}